\documentclass[11pt]{article}
\usepackage{amsfonts, amssymb, latexsym}

\usepackage{color}

\usepackage{color}
\usepackage{mathrsfs}

\def\R{\mathbb{R}}

\def\N{\mathbb{N}}

\def\Z2{{{\mathbb{Z}}^2}}

\def\Z{\mathbb{Z}}
\def\bZ{{\mathbf Z}}

\def\bC{{\mathbf C}}

\def\bV{{\mathbf V}}

\def\S{\mathbb{S}}
\def\0{{\bf 0}}
\def\x{{\bf x}}
\def\y{{\bf y}}
\def\rE{\mathrm{E}}     

\def\v{\vartheta}

\def\R{\mathbb{R}}

\def\N{\mathbb{N}}

\def\Z2{{{\mathbb{Z}}^2}}

\def\Z{\mathbb{Z}}
\def\bZ{{\mathbf Z}}

\def\bC{{\mathbf C}}

\def\S{\mathbb{S}}
\def\0{{\bf 0}}
\def\x{{\bf x}}
\def\y{{\bf y}}
\def\rE{\mathrm{E}}     

\def\v{\vartheta}

\def\R{\mathbb{R}}

\def\N{\mathbb{N}}

\def\Z2{{{\mathbb{Z}}^2}}

\def\Z{\mathbb{Z}}
\def\bZ{{\mathbf Z}}

\def\bC{{\mathbf C}}

\def\bV{{\mathbf V}}

\def\S{\mathbb{S}}
\def\T{\mathbb{T}}
\def\0{{\bf 0}}
\def\x{{\bf x}}
\def\y{{\bf y}}
\def\rE{\mathrm{E}}     

\def\0{{\bf 0}}
\def\rE{\mathrm{E}}     
\def\var{\mathop{\rm var}\nolimits}    
\def\cov{\mathop{\rm cov}\nolimits}    

\usepackage{amsfonts, amssymb, latexsym}
\usepackage{color}
\usepackage{amsfonts, amssymb, latexsym}

\usepackage{epsfig,graphicx}

\def\0{{\bf 0}}
\def\rE{\mathrm{E}}     
\def\var{\mathop{\rm var}\nolimits}    
\def\cov{\mathop{\rm cov}\nolimits}    

\date{}

\begin{document}

\title{Time  Varying   Isotropic   
 Vector  Random  Fields on Spheres }

\author{
   Chunsheng Ma\footnote{  
   Chunsheng  Ma 
   ~~  chunsheng.ma@wichita.edu 
     ~~
Department of Mathematics, Statistics, and Physics, Wichita State
University, Wichita, Kansas 67260-0033, USA;
 School of Economics, Wuhan University of Technology, Hubei 430070,  China;
 ~ and  
 School of Mathematics and Statistics, Hubei Engineering University, 
Xiaogan, Hubei 432000, China 
}
}

\maketitle
\date{}

~ \\

~ \\

\noindent 
{\bf Abstract}    ~~  
For a vector random field that is isotropic and mean square continuous on a sphere and stationary on a temporal domain, this paper derives  a general form of its covariance matrix function and provides   a series representation
for the random field, which  involve the ultraspherical polynomials.  The  series representation  is somehow  an  imitator of the covariance matrix function,  but differs from the
the spectral representation in terms of the ordinary spherical harmonics,  and is    useful for modeling and simulation.
Some semiparametric  models are also illustrated.

~


~

\noindent
{\bf Mathmatics Subject Classification (2010)}   ~~   60G60  $\cdot$ 62M10  $\cdot$  62M30

\newpage

\section{Introduction}

Consider an $m$-variate spatio-temporal  random field $\{ \bZ (\x; t), \x \in \S^d, t \in \T \}$, where $\S^d$  is the spherical shell of radius 1 and center $\0$ in $\R^{d+1}$, i.e.,
$\S^d = \{ \| \x \| =1,  \x \in \R^{d+1} \}$,   $\| \x \|$ is the Euclidean norm of $\x \in \R^{d+1}$,  and $\T$ is either $\R$ or $\Z$.
It is  called a time varying or time dependent  random field on the sphere \cite{Dovidio2014}, \cite{Johns1963}, \cite{Roy1976}.
When $\{ \bZ (\x; t), \x \in \S^d, t \in \T \}$ has finite  second-order  moments,  its mean function and 
   covariance matrix function are given respectively by $\rE \bZ (\x; t)$ and
     $$ \cov (\bZ (\x_1, t_1), \bZ (\x_2;  t_2 )) = \rE \{  (\bZ (\x_1; t_1) - \rE  \bZ (\x_1; t_1))  (\bZ (\x_2; t_2) - \rE  \bZ (\x_2; t_2))' \},  ~~ \x_1, \x_2 \in \S^d, ~ t_1, t_2 \in \T. $$
The primary goal of this paper is to explore the covariance matrix structure and the series representation of an $m$-variate random field       $\{ \bZ (\x; t), \x \in \S^d, t \in \T \}$
that is isotropic and mean square continuous  over the sphere $\S^d$ and stationary over the time domain $\T$, and is mean square continuous on $\S^d \times \T$ if $\T = \R$.

For two points $\x_1$ and $\x_2$ on $\S^d$,
their spherical (angular, or geodesic) distance is 
the distance between  $\x_1$ and $\x_2$ on the largest circle on $\S^d$ that passes
through them; more precisely,
$$ \vartheta (\x_1, \x_2) = \arccos (\x_1' \x_2),  ~~~~ \x_1, \x_2 \in \S^d, $$
or
 $$\vartheta (\x_1, \x_2 ) = \arccos \left( 1- \frac{1}{2} \| \x_1-\x_2 \|^2 \right), ~~~~ \x_1, \x_2 \in  \S^d, $$
 where  $\x'_1\x_2$ is the inner product between $\x_1$ and $\x_2$. Evidently, $0 \le \vartheta (\x_1, \x_2) \le \pi,$ 
  $\S^d$ is a metric space under the  spherical  distance,  and the Euclidean and spherical distances are
closely connected on $\S^d$, with 
        $$  \| \x_1 -\x_2 \| = ( 2 - 2 \x_1' \x_2)^{\frac{1}{2}} = (2- 2 \cos \vartheta (\x_1, \x_2))^{\frac{1}{2}}
            = 2 \sin \left( \frac{\vartheta (\x_1, \x_2)}{2} \right), ~~ \x_1, \x_2 \in \S^d. $$
   
   An $m$-variate random field       $\{ \bZ (\x; t), \x \in \S^d, t \in \T \}$ is said to be (wide-sense)
 isotropic over the sphere $\S^d$ and (wide-sense) stationary over the time domain $\T$, if
 its mean function  $\rE \bZ (\x; t)$ equals a constant vector, and its covariance matrix function $ \cov (\bZ (\x_1;  t_1), \bZ( \x_2;  t_2 ) )$  depends only on the spherical distance $\vartheta (\x_1, \x_2)$ between $\x_1$ and $\x_2$ and
 time lag $t_1-t_2$ between $t_1$ and $t_2$.  In such a case, it covariance matrix function is denoted by $\bC ( \vartheta; t)$, or 
   $$ \bC (\vartheta (\x_1, \x_2); t_1-t_2) = \rE \{  (\bZ (\x_1; t_1) - \rE  \bZ (\x_1; t_1))  (\bZ (\x_2; t_2) - \rE  \bZ (\x_2; t_2))' \} ,   ~~ \x_1, \x_2 \in \S^d, ~ t_1, t_2 \in \T. $$
   It is an $m \times m$ matrix function,  $\bC (\vartheta; -t) = ( \bC (\vartheta; t) )'$, and inequality
    \begin{equation}
     \label{positive.definite}
      \sum_{i=1}^n \sum_{j=1}^n  \mathbf{a}'_i \bC (\vartheta (\x_i, \x_j); t_i-t_j) \mathbf{a}_j \ge 0  
     \end{equation}
      holds for every   $n \in \N$,  any $\x_i \in \S^d$, $t_i \in \T$,  and $\mathbf{a}_i \in \R^m$ ($ i =1, 2, \ldots, n$), where $\N$ stands for the set of positive integers. 
      On the other hand, given an $m \times m$ matrix function with these properties, there exists an $m$-variate Gaussian or elliptically contoured random field   $\{ \bZ (\x; t), \x \in \S^d, t \in \T \}$
      with  $\bC ( \vartheta; t)$ as its covariance matrix function \cite{Ma2011}.
   
   In the scalar case $m=1$,  a spectral analysis is developed in  \cite{Roy1973}, \cite{Roy1976}  when $\T = \Z$ and $d=1$ or $2$,  
   and  a  Fourier series expansion of  $\{ Z(\x;  t), \x \in \S^d, t \in \Z \}$ is derived with  the coefficients being stochastic processes indexed by the time only, as well as   a  spectral representation of its  covariance function.  For $\T = \R$ and $d \ge 2$, the spectral expansion of  a scalar random field 
   $\{ Z(\x;  t), \x \in \S^d, t \in \R \}$  is described by \cite{Mokljacuk1979},
      \begin{equation}
      \label{Mokljacuk1979.exp}
      Z (\x; t) = \sum\limits_{n=0}^\infty \sum_{k=1}^{h(n)}  U_{nk}  (t) S_{n, k}  (\x),  ~~~~  \x \in \S^d,  ~  t \in \R, 
      \end{equation}
     where  $ h (n) = (2n+d-1) \frac{ (n+d-2)!}{(d-1)! n!}$,  $S_{n, k} (\x) $ ($k =1, \ldots, h (n)$) are the orthonormal  spherical harmonics  of degree $n$ on $\S^d$ \cite{Andrews1999}, \cite{Muller1998}, $S_{0, 0} = \frac{1}{\sqrt{ 2 \pi^{ (d+1)/2} \Gamma ((d + 1)/2)}}$,   $\{ U_{nk} (t), t \in \R \}$  ($  n \in  \N_0, k \in \N_0$)  are stationary stochastic processes with
            $$   \cov (U_{n_1 k} (t_1), U_{n_2 j} (t_2)) =  \delta_{n_1 n_2} \delta_{k j} b_{n_1} (t_1-t_2),  ~~~~~~~~ t_1, t_2 \in \T, $$
           $\{ b_n (t),  t \in \R, n \in \N_0 \}$ is a sequence of  stationary covariance functions with 
         $\sum\limits_{n=0}^\infty b_n (0) P_n^{ \left( \frac{d-1}{2} \right)} ( 1) < \infty$,  $\delta_{kj}$ is the  Kronecker symbol, and  $\N_0$ denotes the set of nonnegative integers.
   The covariance function of (\ref{Mokljacuk1979.exp}) is
                 \begin{equation}
      \label{Mokljacuk1979.cov}
       C (\vartheta; t)  =  \sum\limits_{n=0}^\infty  b_n (t) P_n^{ \left( \frac{d-1}{2} \right) } (\cos \vartheta) ,  ~~~~  \vartheta \in [0, \pi],  ~  t \in \R,
      \end{equation}
      where $P_n^{ \left( \frac{d-1}{2} \right) } (x)$ ($n \in \N_0$) are  ultraspherical or   Gegenbauer's  polynomials  \cite{Andrews1999}, \cite{Szego1959}.
         Theoretical investigations and
practical applications of  salar and vector random fields on spheres may be found in   
\cite{Askey1976}, \cite{Bingham1973},  \cite{Cheng2016}, 
\cite{DuMaLi2013}-\cite{Gaspari2006}, \cite{Hannan1970},
 \cite{Johns1963}-\cite{Leonenko2013}, 
           \cite{Ma2012}-\cite{McLeod1986}, 
  \cite{Roy1973}, \cite{Roy1976}, 
    \cite{Yadrenko1983}-\cite{Yaglom1987}.

         For  an $m$-variate random field $\{ \bZ (\x; t), \x  \in \S^d, t \in \T \}$  isotropic and mean square continuous on $\S^d$ and stationary on $\T$,  a general form of its covariance matrix function is  given
         in Section  2,  and a series representation is presented in Section 3.  
         The established  forms  of covariance matrix function and of series representation  mimic each other, and  are  useful for modeling and simulation.  
         Some concluding remarks are made in Section 4, and the theorems are proved in Section 5.

    \section{Covariance  Matrix  Structures}

      For  an  $m$-variate random field  $\{ \bZ (\x; t), \x \in \S^d, t \in \T \}$ that  is 
 isotropic and mean square continuous  over  $\S^d$ and stationary on $\T$,
  its covariance matrix function $\bC (\vartheta; t)$  is a continuous  function of $\vartheta \in [0, \pi]$, and is also a continuous function of $t \in \R$ if $\T = \R$.
 This section derives the general form of    $\bC (\vartheta; t)$, which involves  ultraspherical polynomials.
 
  We  start with a brief review of  ultraspherical polynomials,  and refer to \cite{Andrews1999} and \cite{Szego1959} for comprehensive treatments.   
     For $\lambda > 0$,
  the  ultraspherical or Gegenbauer's polynomials, $P^{(\lambda)}_n (x)$,
   $n \in \N_0$,  are the coefficients of $u^n$ in the power series expansion of the function 
  $(1-2 u x +u^2)^{-\lambda}$, i.e.,
         \begin{equation}
         \label{Gpolynomial}
          (1-2 u x +u^2)^{-\lambda} = \sum_{n=0}^\infty u^n P^{(\lambda)}_n (x), ~~~~~  x \in \R,  ~  | u |  < 1.  
          \end{equation}   
        They     can be alternatively defined  through the recurrence formula
             \begin{eqnarray*}
             \left\{
             \begin{array}{lll}
           P^{(\lambda)}_0 (x) & \equiv & 1,   \\
           P^{(\lambda)}_1 (x)  & = &    2 \lambda x,           \\
           P^{(\lambda)}_n (x)  &  =   &   \frac{2 ( \lambda+n-1) x  P^{(\lambda)}_{n-1} (x) - (2 \lambda+n-2) P^{(\lambda)}_{n-2} (x)}{
                                                           n},   ~~~~~ x \in \R,  ~  n \ge 2.
              \end{array}    \right.                                               
           \end{eqnarray*}
           Some special cases and particular values are 
         \begin{eqnarray*}
             &   & P^{(1)}_n ( \cos \vartheta)   =   \frac{\sin ( (n+1) \vartheta )}{\sin \vartheta},  ~~~~~~ \vartheta \in [0, \pi],    \\
            &  & P^{(\lambda)}_n (1)   =     { 2 \lambda+n-1 \choose n},
           \end{eqnarray*}
           and
           $$  \left| P^{(\lambda)}_n (x)  \right| \le P^{(\lambda)}_n (1),   ~~~~~~ |x| \le 1. $$
           In the particular case  $\lambda =\frac{1}{2}$, $P^{ \left( \frac{1}{2} \right)}_n (x)$  ($n \in \N_0$) are
          the Legendre polynomials. 
          
   The ultraspherical polynomials  are polynomial solutions of  the differential equation
        $$ (1-x^2) \frac{d^2 y}{d x^2} - (2 \lambda+1) x \frac{d y}{d x} +n (2 \lambda+n) y  = 0, $$
        and possess two types of  orthogonal  properties. First, they 
         are orthogonal with respective to the weight function $(1-x^2)^{\lambda-\frac{1}{2}}$ on $[-1, 1]$, in the sense that
           \begin{equation}
            \label{Gegenbauer.orth}
             \int_{-1}^1 P^{(\lambda)}_i (x) P^{(\lambda)}_j (x)  (1-x^2)^{\lambda-\frac{1}{2}} d x 
                =  \left\{
                \begin{array}{ll}
                \frac{\pi 2^{1- 2 \lambda} \Gamma (i+2 \lambda)}{i! (\lambda+i) (\Gamma (\lambda))^2},
                 ~   &  ~ i =j, \\
                 0, ~ & ~ i \neq j.
                 \end{array}   \right. 
                \end{equation}   
              Second, they are  orthogonal  over     $\S^d$ ($d \ge 2$), as the following lemma describes, which 
              is a special case of the   Funk-Hecke formula  (\cite{Andrews1999}, \cite{Muller1998})  that is useful in simplifying calculations of certain integrals over  $\S^d$.
              
              ~
                 
                   \noindent
 {\bf  Lemma 1}  {\em For  $i, j \in \N_0$, if $d \ge 2$, then }
            \begin{eqnarray*}
             \int_{\S^d}  P_i^{\left( \frac{d-1}{2} \right)} (\x'_1 {\bf z})  P_j^{\left( \frac{d-1}{2} \right)} (\x'_2 {\bf z})d {\bf z} =
             \left\{
             \begin{array}{ll}
                \frac{  \omega_d}{\alpha_i^2}  P_i^{ \left( \frac{d-1}{2}  \right)} (\x'_1 \x_2 ) ,   ~   &  ~ i = j, \\
              0, ~   &  ~ i \neq j,
              \end{array}   \right. 
              \end{eqnarray*}
      {\em         where  $\omega_{d} =\frac{2 \pi^{\frac{d+1}{2}}}{\Gamma \left( \frac{d+1}{2} \right)}$ is the surface area of  $\S^{d}$, and   }
                  \begin{equation}
         \label{alpha.n}
          \alpha_n  =  \left(  \frac{2n+d-1}{d-1} \right)^{\frac{1}{2}},    ~~~~~~   n \in  \N_0.
          \end{equation}

   ~

            In terms of  the orthonormal spherical harmonics,  $P_n^{ \left( \frac{d-1}{2} \right) } (\cos \vartheta)$ or $P_n^{ \left( \frac{d-1}{2} \right) } (\x' \y)$ can be expressed as 
         (see, e.g.,  Theorem 9.6.1 of \cite{Andrews1999})
            \begin{equation}
        \label{Andrews1999}
        P_n^{ \left( \frac{d-1}{2} \right)}  (\x' \y ) 
          =  \frac{\omega_d}{\alpha_n^2}    \sum_{k=1}^{h(n)}   S_{n, k}  (\x)    S_{n, k}  (\y),
           ~~~~~ \x, \y \in \S^d,  ~ n \in \N,
         \end{equation} 
         from which its positive definiteness follows directly, noticing that $\frac{h(n)}{P_n^{ \left( \frac{d-1}{2} \right)}  (1)  } = \alpha_n^2$.
         The elementary positive definite spherical functions on $\S^d$   are the positive scalar products  \cite{Cartan1929}  of 
          $P_n^{ \left( \frac{d-1}{2} \right) } (\cos \vartheta)$, $ n \in \N_0$, which actually  form a basis  \cite{Schoenberg1942} of the set of isotropic, continuous, and positive definite functions on $\S^d$.
           A probability interpretation 
          for these   elementary positive definite functions on the sphere is provided in Lemma 2 below \cite{Ma2016},  which  illustrates  a basis of the set of isotropic and mean square random fields on $\S^d$. It  will be employed in the proofs of Theorems 1, 2 and 4. 
 
          ~

\noindent
{\bf Lemma 2}   {\em   If
${\bf U}$ is a 
       (d+1)-dimensional random vector uniformly distributed on $\S^d$ ($ d \ge 2$),  then, for a fixed  }   $n \in \N$,
        \begin{equation}
        \label{ele.rf}
             Z_n (\x) = \alpha_n P_n^{ \left( \frac{d-1}{2} \right) } (\x' {\bf U}),  ~~~~ \x \in \S^d,  
        \end{equation}     
     {\em   is an isotropic random field  with mean  
       0 and covariance function   }
       \begin{equation}
       \label{ultra.cov}
        \cov  (  Z_n (\x_1),   Z_n (\x_2 ) )
                =  P_n^{ \left( \frac{d-1}{2} \right) } ( \cos \vartheta (\x_1, \x_2)),
                       ~~~ \x_1, \x_2 \in \S^d, 
       \end{equation}                
     {\em    where $\alpha_n$ is defined in (\ref{alpha.n}).
         Moreover,  for $i \neq j$, 
       $\{ Z_i (\x), \x \in \S^d \}$ and  $\{ Z_j (\x), \x \in \S^d \}$ are uncorrelated; that is  }
       $$ \cov ( Z_i (\x_1), Z_j (\x_2)) =0,   ~~~~~~~  \x_1, \x_2 \in  \S^d. $$

 ~
 
 Alternatively, assume that $Z_1, \ldots, Z_{h(n)}$ are uncorrelated random variables with mean 0 and variance
  1. Then
      \begin{equation}
      \label{stoc.ultrasph}
       Z_n (\x) = \frac{\omega_d^{\frac{1}{2}}}{\alpha_n}  \sum_{k=1}^{h(n)} Z_k   S_{n, k} (\x),  ~~~~~~~ \x \in \S^d,
       \end{equation}
       is an isotropic random field with mean 0 and covariance function $P_n^{ \left( \frac{d-1}{2} \right) } ( \cos \vartheta (\x_1, \x_2))$; see
    page 77 of \cite{Yadrenko1983}. 
    More interestingly, (\ref{ele.rf}) may be thought of as a special case of (\ref{stoc.ultrasph}) by selecting
     $$ Z_k = \omega_d^{\frac{1}{2}}  S_{n, k} (\mathbf{U}),   ~~~ k = 1, \ldots, h(n), $$
     with the help of identity (\ref{Andrews1999}).
          
 ~

    \noindent
    {\bf Theorem 1}  ~~   {\em     If  an  $m$-variate random field  $\{ \bZ (\x; t), \x \in \S^d, t \in \T \}$ is     isotropic and mean square continuous  over  $\S^d$ and stationary on $\T$, then 
     $ \frac{\bC (\vartheta; t) + \bC (\vartheta; -t)}{2} $ is of the form   }
       \begin{equation}
       \label{cov.matrix.fun.1}
        \frac{\bC (\vartheta; t) + \bC (\vartheta; -t)}{2}   =
         \left\{
          \begin{array}{ll}
          \sum\limits_{n=0}^\infty  \mathbf{B}_n (t) \cos (n  \vartheta), ~    &   ~  d =1,   \\
          ~   &  ~ ~~~~~~~~~   \vartheta \in [0, \pi], ~ t \in \T,    \\
          \sum\limits_{n=0}^\infty  \mathbf{B}_n (t) P_n^{ \left( \frac{d-1}{2} \right) } (\cos \vartheta),    ~ &   ~  d \ge 2,   
          \end{array}
          \right.
       \end{equation}
   {\em   where, for each fixed $t \in \T$,   $ \mathbf{B}_n (t)$ ($ n \in \N_0$) are $m \times m $ symmetric matrices and  $\sum\limits_{n=0}^\infty  \mathbf{B}_n (t)$ ($d =1$) or $\sum\limits_{n=0}^\infty  \mathbf{B}_n (t)  P_n^{ \left( \frac{d-1}{2} \right) } (1)$ ($d \ge 2$) converges,  
     and, for each fixed $n \in \N_0$,  $ \mathbf{B}_n (t)$ is a stationary covariance matrix function on $\T$.  }
     
 {\em    In particular,  when $\bC (\vartheta; t)$ is spatio-temporal symmetric in the sense  that  }
               $$ \bC ( \vartheta;  - t  ) =\bC ( \vartheta;  t  ),   ~~~~~~~~ \vartheta \in [0, \pi], ~ t \in \T, $$
      {\em       it  takes  the form } 
     $$   \bC (\vartheta; t)    =
         \left\{
          \begin{array}{ll}
          \sum\limits_{n=0}^\infty  \mathbf{B}_n (t) \cos (n  \vartheta), ~    &   ~  d =1,   \\
          ~   &  ~ ~~~~~~~~~   \vartheta \in [0, \pi], ~ t \in \T,    \\
          \sum\limits_{n=0}^\infty  \mathbf{B}_n (t) P_n^{ \left( \frac{d-1}{2} \right) } (\cos \vartheta),    ~ &   ~  d \ge 2.  
          \end{array}
          \right.  $$

       ~

       In the next theorem $m \times m$ matrices  $ \mathbf{B}_n (t)$ ($ n \in \N_0$) are not necessarily   symmetric. 
       One simple such  example is 
            $$  \mathbf{B} (t) = \left\{
                                           \begin{array}{ll}
                                           \mathbf{\Sigma}+\mathbf{\Phi}        \mathbf{\Sigma}  \mathbf{\Phi}',  ~  &  ~ t =0, \\
                                            \mathbf{\Phi}        \mathbf{\Sigma},   ~   &  ~  t = -1, \\
                                            \mathbf{\Sigma}  \mathbf{\Phi}',  ~  &  ~ t = 1, \\
                                            \mathbf{0},  ~ & ~  t = \pm 2, \pm 3, \ldots,
                                            \end{array}    \right. $$
    which is     the covariance matrix function of  an $m$-variate first order moving average time series $\bZ (t) = \mathbf{\varepsilon}  (t) + \mathbf{\Phi}  \mathbf{\varepsilon} (t-1),
       t \in \Z, $  where $\{  \mathbf{\varepsilon}  (t), t \in \Z \}$ is $m$-variate  white noise with $\rE \mathbf{\varepsilon}  (t) = \mathbf{0}$ and
       $\var ( \mathbf{\varepsilon}  (t)) = \mathbf{\Sigma}$, and $\mathbf{\Phi}$ is an $m \times m$ matrix. 
       
       ~
       
        \noindent
    {\bf Theorem 2}   ~~ (i) {\em An  $m \times m$ matrix function}  
         \begin{equation}
             \label{cov.matrix.fun.2.0}
        \bC (\vartheta; t)    =
          \sum\limits_{n=0}^\infty  \mathbf{B}_n (t) \cos (n  \vartheta), ~   ~~   \vartheta \in [0, \pi], ~ t \in \T,   
          \end{equation}
       {\em is the covariance matrix function of   an  $m$-variate Gaussian or elliptically contoured random field
       $\{ \bZ (\x; t), \x \in \S^1, t \in \T \}$  if and only if $\sum\limits_{n=0}^\infty  \mathbf{B}_n (0)$   converges and 
        $ \mathbf{B}_n (t)$ is a stationary covariance matrix function on $\T$ for each fixed $n \in \N_0$.} 
        
        (ii)  {\em Let $d \ge 2$.  An  $m \times m$ matrix function}  
         \begin{equation}
             \label{cov.matrix.fun.2.1}
        \bC (\vartheta; t)    =
          \sum\limits_{n=0}^\infty  \mathbf{B}_n (t)  P_n^{ \left( \frac{d-1}{2} \right) } (\cos \vartheta), ~   ~~   \vartheta \in [0, \pi], ~ t \in \T,   
          \end{equation}
       {\em is the covariance matrix function of   an  $m$-variate Gaussian or elliptically contoured random field
        on $ \S^d \times  \T $  if and only if  $\sum\limits_{n=0}^\infty  \mathbf{B}_n (0) P_n^{ \left( \frac{d-1}{2} \right) } (1)$   converges and
        $ \mathbf{B}_n (t)$ is a stationary covariance matrix function on $\T$ for each fixed $n \in \N_0$.} 
       
       ~

            Gaussian and second-order elliptically contoured  random fields  form one of the largest sets, if not the largest set, 
          which allow any possible correlation structure \cite{Ma2011}. 
          The covariance matrix functions developed in Theorem 2  can be adopted
          for a Gaussian or elliptically contoured vector random field.
          However, they may not be available for other non-Gaussian  random fields, such as  a log-Gaussian, $\chi^2$, K-distributed, or skew-Gaussian one,  for which  admissible correlation structure must be investigated on a case-by-case basis.

          ~

\noindent
{\bf Example 1}  Given  an $m \times m$  matrix function $\mathbf{B} (t), t \in \T,$ with all entries $b_{ij} (t)$ less than 1   in absolute value, consider an $m \times m$ matrix function
$\bC (\vartheta; t)$ with entries
   $$ C_{ij} (\vartheta; t) =  -\ln \left\{  \frac{1}{2} \left[  1-b_{ij} (t) \cos \vartheta+ \left( 1-2 b_{ij} (t) \cos \vartheta+b^2_{ij} (t)  \right)^{\frac{1}{2}}  \right]  \right\},  ~~~~ \vartheta \in [0, \pi], 
         ~ t \in \T,   $$
         \hfill   $~ i, j =1, \ldots, m.$   ~~~~
         
         \noindent
         It is  the covariance matrix function of   an  $m$-variate Gaussian or elliptically contoured random field  $\{ \bZ (\x; t), \x \in \S^2, t \in \T \}$ if and only if 
         $ \mathbf{B} (t)$ is a stationary covariance matrix function on $\T$.
         In fact,  a  version  (\ref{cov.matrix.fun.2.1})   of $\bC(\vartheta; t)$ can be established by taking  $ \mathbf{B}_0 (t) \equiv  \mathbf{0}, $ $ \mathbf{B}_n (t) =\frac{1}{n} (\mathbf{B} (t))^{\circ n}, n \in \N,$ and using the identity
   (see, e.g.,  (5) on page 128 of \cite{Mangulis1965})
   $$ \sum_{n=1}^\infty \frac{u^{n}}{n }  P_n^{\left( \frac{1}{2} \right)} ( x) =
   -\ln \left\{ \frac{1}{2} \left[   1-u x+(1-2 u x+u^2)^{\frac{1}{2}}  \right] \right\}, ~~~~~~~~ |x| \le 1, ~ |u| < 1, $$
    where  ${\bf B}^{\circ p}$ denotes the  Hadamard  $p$ power of ${\bf B} =(b_{ij})$, whose entries are $b_{ij}^p$, the $p$ power of $b_{ij}, i, j = 1, \ldots, m$.

       ~

       A covariance matrix function $\bC (\vartheta; t)$ defined on $\S^d \times \T$ is also a covariance matrix function on $\S^{d_0} \times \T$,  provided that $1 \le  d_0 < d$,  just as a point 
       $\x \in  \S^{d_0}$ can be thought of as a point $(\x, \mathbf{0})' $ on $\S^d$. 
       A  covariance matrix function $\bC (\vartheta; t)$ on all $\S^d \times \T$ ($ d \in \N$) is called a covariance matrix function on $\S^\infty \times \T$,
        with $\S^\infty$  being an infinite dimensional sphere in Hilbert space.
       A general form of this type of covariance matrix structures is given next. 
       
       ~

    \noindent
    {\bf Theorem 3}  ~~  (i)  {\em   If  an  $m$-variate random field  $\{ \bZ (\x; t), \x \in \S^\infty, t \in \T \}$ is     isotropic and mean square continuous  over  $\S^\infty$ and stationary on $\T$, then 
      }
       \begin{equation}
       \label{cov.matrix.fun.3}
        \frac{\bC (\vartheta; t) + \bC (\vartheta; -t)}{2}   =
            \sum\limits_{n=0}^\infty  \mathbf{B}_n (t) \cos^n  \vartheta, 
            ~~~~   \vartheta \in [0, \pi], ~ t \in \T,   
       \end{equation}
  {\em    where, for each fixed $t \in \T$,   $ \mathbf{B}_n (t)$ ($ n \in \N_0$) are $m \times m $ symmetric matrices and  $\sum\limits_{n=0}^\infty  \mathbf{B}_n (t)  $ converges,  
     and, for each fixed $n \in \N_0$,  $ \mathbf{B}_n (t)$ is a stationary covariance matrix function on $\T$. }
     
  {\em    In particular,   a   spatio-temporal symmetric  $\bC (\vartheta; t)$    is of the form   }
             $$ \bC (\vartheta; t)    =
          \sum\limits_{n=0}^\infty  \mathbf{B}_n (t) \cos^n  \vartheta, ~   ~~~~   \vartheta \in [0, \pi], ~ t \in \T.   $$
       
  (ii)  {\em     An $m \times m$ matrix function 
               \begin{equation}
             \label{cov.matrix.fun.4}
        \bC (\vartheta; t)    =
          \sum\limits_{n=0}^\infty  \mathbf{B}_n (t) \cos^n  \vartheta, ~   ~~~~   \vartheta \in [0, \pi], ~ t \in \T,
       \end{equation}
       is the covariance matrix function of 
   an $m$-variate Gaussian or elliptically contoured random field
       on $ \S^\infty \times  \T $  if and only if   $\sum\limits_{n=0}^\infty  \mathbf{B}_n (0)  $ converges and    $ \mathbf{B}_n (t)$ is a stationary covariance matrix function on $\T$ for each fixed $n \in \N_0$. }

       ~
       
     One may use Lemma 1 of  \cite{Schoenberg1942}  to deduce (\ref{cov.matrix.fun.3}).  Instead,   a more efficient  approach   is based on  Lemma 3 below,  
      which expresses  $x^n$  as a convex combination of ultraspherical polynomials,   explains a close connection between $\cos^n \vartheta$ and $P_{k}^{ \left( \frac{d-1}{2} \right) } ( \cos \vartheta)$
      ($k =0, 1, \ldots, n$),
      where the former is for the basis of  the covariance matrix structure on $\S^\infty$ and the latter on $\S^d$, and provides an approach to generate an isotropic random field on $\S^d$ with covariance function
      $\cos^n \vartheta$. 
       
       ~
       
       \noindent
       {\bf Lemma 3}    {\em  Let $n \in \N_0$. } 
       
         \begin{itemize}
         \item[(i)]    $x^n$ {\em can be expressed as  }
        \begin{equation}
        \label{Binghamlemma1973}
        x^n = \sum\limits_{k=0}^{ \left[ \frac{n}{2} \right]}
                   \beta_{k,n}^{\left( \frac{d-1}{2} \right)} P_{n-2 k}^{ \left( \frac{d-1}{2} \right) } (x),    ~~~~~~~ |x| \le 1, 
         \end{equation}
     {\em     where $[u]$ denotes the  integer part of a real number $u$, and }
           $$        \beta_{k,n}^{\left( \frac{d-1}{2} \right)}  =  \frac{n! \left( n-2k +\frac{d-1}{2} \right) \Gamma \left( \frac{d-1}{2} \right) }{
                                      2^n k!  \Gamma \left( n-k +\frac{d+1}{2} \right)  },   ~~~~~~~~~~~~  k = 0, 1, \ldots, \left[ \frac{n}{2} \right];  $$   
     
         \item[(ii)]   $\cos^n \vartheta$  {\em can be expressed as  }
        \begin{equation}
        \label{Binghamlemma1}
        \cos^n \vartheta  = \sum\limits_{k=0}^{ \left[ \frac{n}{2} \right]}
                   \beta_{k,n}^{\left( \frac{d-1}{2} \right)} P_{n-2 k}^{ \left( \frac{d-1}{2} \right) } (\cos \vartheta ),    ~~~~~~~    \vartheta \in [0, \pi];
         \end{equation}
         
         \item[(iii)]    {\em   If
${\bf U}$ is a 
       (d+1)-dimensional random vector uniformly distributed on $\S^d$ ($ d \ge 2$),  then }
        \begin{equation}
        \label{Binghamlemma2}
             Z (\x) =    \sum\limits_{k=0}^{ \left[ \frac{n}{2} \right]}  \alpha_k 
                  \left(  \beta_{k,n}^{\left( \frac{d-1}{2} \right)} \right)^{\frac{1}{2} }   P_{n-2 k}^{ \left( \frac{d-1}{2} \right) } (\x' {\bf U}),  ~~~~ \x \in \S^d,  
        \end{equation}     
     {\em   is an isotropic random field  with mean  
       0 and covariance function   $\cos^n \vartheta$. }

         \end{itemize}
         
         Identity  (\ref{Binghamlemma1973})  is an alternative  version of  Lemma 1  of  \cite{Bingham1973}, and Part (iii) of Lemma 3 follows from Lemma 2 and (\ref{Binghamlemma1}).  Another method  generating  an isotropic random field  on $\S^d$  with covariance function   $\cos^n \vartheta$ is presented in         
   Subsection 12.3 of  \cite{Cohen2012}. 
   
   ~
   
   \noindent
   {\bf Example 2}  
 An $m \times m$ matrix function $\bC (\vartheta; t)$ whose entries are  second order polynomials of $\vartheta$, 
          \begin{equation}
           \label{ex1}
             \bC(\vartheta; t) = {\bf B}_0 (t) + {\bf B}_1 (t) \vartheta +{\bf B}_2 (t)  \vartheta^2,  ~~~~~~~~~~ \vartheta  \in [0, \pi],  ~ t \in \T, 
          \end{equation}   
                     is a covariance matrix function on $\S^\infty \times \T$ if and only if 
                   ${\bf B}_2(t) $,  $-{\bf B}_1 (t)- \pi {\bf B}_2 (t) $, and   ${\bf B}_0 (t) +\frac{\pi}{2} {\bf  B}_1 (t)  +\frac{\pi^2}{4} {\bf B}_2(t)$
                   are   stationary covariance matrix functions on $\T$.
 To  apply Theorem 3 to the function (\ref{ex1}), we  employ  the formula 
      $$ \vartheta = \frac{\pi}{2} - \arcsin (\cos \vartheta),   ~~~~~~~ \vartheta \in [0, \pi],  $$
and              
  the Taylor  expansions of $\arcsin x$ and $(\arcsin x)^2$,
     $$ \arcsin x = \sum\limits_{n=0}^\infty \frac{(2n)!}{2^{2n} (n!)^2 (2n+1)} x^{2n+1},  ~~~~~~~~~ | x | \le 1,   $$ 
     $$   (\arcsin x)^2= \sum_{n=1}^\infty  \frac{2^{2n-1} ( (n-1)!)^2}{(2n)!}  x^{2n}, ~~~~~~~~~~  | x | \le 1,      $$ 
   and obtain  a version (\ref{cov.matrix.fun.4}) of $\bC (\vartheta; t)$,
         \begin{eqnarray*}
         &   & \bC  (\vartheta; t)  =   {\bf B}_0 (t) + {\bf B}_1 (t)  \left( \frac{\pi}{2} - \arcsin (\cos \vartheta) \right)  +{\bf B}_2 (t)  \left( \frac{\pi}{2} - \arcsin (\cos \vartheta)   \right)^2 \\
                                 & = &    {\bf B}_0 (t)  +    \frac{\pi}{2}    {\bf B}_1  (t)    +  \frac{\pi^2}{4}  {\bf B}_2 (t)
                                            - \left( {\bf B}_1 (t) + \pi {\bf B}_2  (t)  \right) \arcsin  (\cos \vartheta) 
                                             + {\bf B}_2  (t)  (\arcsin   (\cos \vartheta)   )^2 \\
                                  & = &      {\bf B}_0 (t) +    \frac{\pi}{2}    {\bf B}_1  (t)   +  \frac{\pi^2}{4}  {\bf B}_2 (t)
                                            - \left( {\bf B}_1 (t) + \pi  {\bf B}_2  (t)  \right)   \sum\limits_{n=0}^\infty \frac{(2n)!}{2^{2n} (n!)^2 (2n+1)}  \cos^{2n+1} \vartheta  \\
                                  &   &           + {\bf B}_2 (t)  \sum_{n=1}^\infty  \frac{2^{2n-1} ( (n-1)!)^2}{(2n)!}   \cos^{2n} \vartheta,  
                                     ~~ ~ \vartheta  \in [0, \pi], ~ t \in \T, 
         \end{eqnarray*}
         whose coefficients are  stationary covariance matrix functions on $\T$   if and only if 
            ${\bf B}_2(t) $,  $-{\bf B}_1 (t)- \pi {\bf B}_2 (t) $, and   ${\bf B}_0 (t) +\frac{\pi}{2} {\bf  B}_1 (t)  +\frac{\pi^2}{4} {\bf B}_2(t)$
                   are so. 
                   
                   In particular, in (\ref{ex1}) taking ${\bf B_2} (t)  \equiv {\bf 0}$ yields that
                  $$  \bC(\vartheta;  t) = {\bf B}_0 (t)  + {\bf B}_1 (t)  \vartheta,  ~~~~~~~~~~ \vartheta  \in [0, \pi], ~  t \in \T, $$
                  is  a covariance matrix function on $\S^\infty \times \T$ 
if and only if ${\bf B}_0 (t) + \frac{\pi}{2}  {\bf B}_1 (t) $  and $-{\bf B}_1 (t) $ are     stationary covariance matrix functions on $\T$.
 Moreover, under these conditions, 
         $$ \bC (\vartheta; t) = ( {\bf B}_0 (t) -  {\bf B}_1 (t) \vartheta )^{\circ p},  ~~~~~~~~~~ \vartheta  \in [0, \pi],  ~ t \in \T,  $$
      is  also a  covariance matrix function by Theorem 6 of \cite{Ma2011},      where $p $ is a natural number.  
      
      ~

      \noindent
      {\bf Example 3}  An $m \times m$ matrix function $\bC (\vartheta; t) $ with entries
       \begin{equation}
       \label{ex2}
       C_{ij} (\vartheta; t) = \exp \left(  -\frac{\pi}{2} b_{ij} (t)   - b_{ij} (t)  \vartheta  \right),    ~~~~~~~ \vartheta \in [0, \pi],  ~ t \in \T, ~~  i,  j =1, \ldots, m,   
       \end{equation}
   is a  covariance matrix function  on $\S^\infty \times \T$ if and only if     ${\bf B}(t) $
   is a  stationary covariance matrix function on $\T$.
 Theorem 3 is applicable, after  we use the Taylor series  of $\exp (\arcsin x)$ (see, for instance, formula 1.216 of  \cite{Gradshteyn2007}),
    $$ \exp ( \arcsin  x) =  \sum_{n=0}^\infty \frac{(n+1) x^n}{n!},        ~~~~~~~~~~~~~~~~ |x| \le 1,  $$
    to represent (\ref{ex2}) as the form of (\ref{cov.matrix.fun.4}),
       \begin{eqnarray*}
        \bC      ( \vartheta; t )  &   =  &  \bC \left(    \frac{\pi}{2} - \arcsin (\cos \vartheta); t     \right)    \\ 
                  &    = &   \sum_{n=0}^\infty \frac{n+1}{n!} {\bf B}^{\circ n} (t)  \cos^n \vartheta,     ~~~~~~~~~~~~~~~~ \vartheta \in [0, \pi], ~  t \in \T,  
        \end{eqnarray*}              
   whose coefficients are stationary covariance matrix functions on $\T$  if and only if ${\bf B} (t)$ is so.

\section{Series  Representations }

 For  an $m$-variate random field  with covariance matrix function (\ref{cov.matrix.fun.2.0}) or  (\ref{cov.matrix.fun.2.1})  this section 
 provides a series representation,  which is a mimic of  (\ref{cov.matrix.fun.2.0}) or (\ref{cov.matrix.fun.2.1})  involving ultraspherical polynomials. 
 A purely spherical version is given in \cite{Ma2016}. 
 Two cases  $d \ge 2$ and $d=1$ are  treated  in Theorems 4 and 5 separately, since the main tool for the construction,   Lemma 2,  applies to the case $d \ge 2$ only. 
 
 ~

\noindent
{\bf Theorem 4}
{\em  Assume that 
 $\{ \bV_n (t),  t \in \T   \}$ is an $m$-variate  stationary stochastic process   with 
$\rE  \bV_n= {\bf 0}$ and  $\cov ( \bV_n (t_1), \bV_n (t_2) ) =  \alpha_n^2 \mathbf{B}_n (t_1-t_2)$ for each fixed $n \in \N_0$,     ${\bf U}$  is a 
       (d+1)-dimensional random vector uniformly distributed on $\S^d$ ($ d \ge 2$),  and $\mathbf{U}$ and   $\{ \bV_n (t),  t \in \T   \}$,  $n \in \N_0$, are independent.
      If  $\sum\limits_{n=0}^\infty {\bf B}_n (0)  P_n^{ \left( \frac{d-1}{2} \right) } (1)$ converges,  then an  $m$-variate  random field  }
       \begin{equation}
     \label{stoc1}
      \bZ (\x; t) = \sum_{n=0}^\infty   \bV_n (t)  P_n^{ \left( \frac{d-1}{2} \right) } (\x' {\bf U}),
      ~~~~~~ \x   \in \S^d, ~ t \in \T,
      \end{equation}
    {\em   is isotropic and mean square continuous on $\S^d$,  stationary on $\T$,  and possesses   mean   ${\bf 0 }$ and  covariance matrix function  }
 (\ref{cov.matrix.fun.2.1}).

~

        The distinct terms of (\ref{stoc1}) are uncorrelated each other, according to Lemma 2 and the independent assumption among $\mathbf{U}, \mathbf{V}_i (t), \mathbf{V}_j (t)$,
           $$ \cov  \left(   \bV_i (t) P_i^{ \left( \frac{d-1}{2} \right) } (\x' {\bf U}),  ~
                          \bV_j  (t) P_j^{ \left( \frac{d-1}{2} \right) } (\x' {\bf U})  \right) = \0, ~~~ \x \in \S^d, ~ t \in \T,  i \neq j. $$
 To see  how   $\bZ (\x; t)$, ${\bf V}_n (t)$ and ${\bf U}$ are related to each other in (\ref{stoc1}),  we multiply both sides of (\ref{stoc1}) by $P_n^{ \left( \frac{d-1}{2} \right) } (\x' {\bf U})$,   integrate  over $\S^d$, and  obtain, in view of Lemma 2, 
     \begin{eqnarray*}
     &       &  \int_{\S^d}  \bZ (\x; t) P_n^{ \left( \frac{d-1}{2} \right) } (\x' {\bf U}) d \x  \\
     &  =  &   \sum_{k=0}^\infty  \bV_n  (t)  \int_{\S^d} P_k^{ \left( \frac{d-1}{2} \right) } (\x' {\bf U})
                       P_n^{ \left( \frac{d-1}{2} \right) } (\x' {\bf U}) d \x  \\
     &  =  &  \frac{\omega_d}{\alpha_n^2}     P_n^{ \left( \frac{d-1}{2} \right) } (1) \bV_n (t)  ,
      \end{eqnarray*}
   or
       $$  \bV_n (t) =   \frac{\alpha^2_n}{ \omega_d P_n^{ \left( \frac{d-1}{2} \right) } (1)}   
             \int_{\S^d}  \bZ (\x; t) P_n^{ \left( \frac{d-1}{2} \right) } (\x' {\bf U}) d \x,  ~~~~~ t \in \T, ~  n \in \N_0. $$
     
     ~
     
     \noindent
{\bf Example 4} 
 Suppose that $\mathbf{B} (t)$ is an $m \times m$  stationary covariance matrix function on $\T$ with all entries less than 1 in absolute value,  all entries of $\mathbf{B}_0 (t)$ equal 1, and $\mathbf{B}_n (t) =   (\mathbf{B} (t))^{\circ n}, n \in \N$.
Then (\ref{stoc1}) defines  an $m$-variate isotropic random field on $\S^d \times \R$ ($d \ge 2$),  with  mean   ${\bf 0 }$ and   direct/cross covariance  functions
     $$ C_{ij} ( \vartheta; t) =  \left( 1-2 b_{ij} (t) \cos \vartheta +b_{ij}^2 (t) \right)^{-\frac{d-1}{2}},  ~~~~ \vartheta \in [0, \pi], ~ t \in \T, ~ ~ i, j =1, \ldots, m, $$
which follows from (\ref{cov.matrix.fun.2.1}) and  the expansion (\ref{Gpolynomial}).
     
~

Theorem 4  does not apply to the unit circle   case $d=1$,  just as Lemma 2 is limited to $d \ge 2$.
 To deal with the case $d =1$,   for two points $\x_1$ and  $\x_2$ on the unit circle $\S^1$,  denote their Cartesian coordinates by
    $\x_k = ( \cos \theta_k, \sin \theta_k)'$, respectively, where $\theta_k$ is the angular coordinate of $\x_k$ in polar coordinates with $0 \le \theta_k \le 2 \pi,$  $k=1, 2$.
    In terms of $\theta_1$ and $\theta_2$,  the angular distance $\vartheta (\x_1, \x_2)$ between $\x_1$ and $\x_2$  can be expressed as
    $$   \vartheta (\x_1, \x_2) = \min ( | \theta_1-\theta_2|, ~ 2 \pi - | \theta_1-\theta_2| ),  $$
  since $\x'_1 \x_2 = \cos (\vartheta (\x_1, \x_2)) = \cos (\theta_1-\theta_2)$.
  The following theorem provides a series  representation for an $m$-variate random field that is isotropic  and mean square continuous  on the unit circle 
  and stationary on $\T$. 
  
      ~

\noindent
{\bf Theorem 5}
 {\em Suppose that, for each $n \in \N_0$,     $\{ \bV_{n1} (t),  t \in \T \}$ and $\{ \bV_{n2} (t), t \in \T \}$ are   $m$-variate  stationary  stochastic processes with  
$\rE  \bV_{nk} (t) = {\bf 0}$ and  $\cov ( \bV_{nk} (t_1), \bV_{nk} (t_2) ) =  \mathbf{B}_n (t_1-t_2)$, $k=1, 2$,
and that  
 $\{ \bV_{n1} (t), t \in \T \}$ and  $\{ \bV_{n2} (t),  t \in \T \}$ ($  n  \in \N_0 $) are independent.
  If $ \sum\limits_{n=0}^\infty  {\bf B}_n (0)$ converges,  then  }
       \begin{equation}
     \label{stoc2}
      \bZ (\x; t ) =\sum_{n=0}^\infty  (  \bV_{n1} (t) \cos (n \theta)  + \bV_{n2}  (t)  \sin (n \theta)   ),  ~~~
          \x = (\cos \theta, \sin \theta)' \in \S^1,  ~  t \in \T, 
      \end{equation}
   {\em     is an $m$-variate  random field on $\S^1 \times \T$,  with  mean   ${\bf 0 }$ and  covariance matrix function  }
 (\ref{cov.matrix.fun.2.0}).  
 
 ~
 
     \noindent
{\bf Example 5}  Let $\mathbf{B} (t)$ be an $m \times m$  stationary covariance matrix function on $\T$, all entries of $\mathbf{B}_0 (t)$ equal 1, and $\mathbf{B}_n (t) =  \frac{1}{n!} (\mathbf{B} (t))^{\circ n}, n \in \N$.
Then (\ref{stoc2}) is  an $m$-variate isotropic random field on $\S^1 \times \R$,  with  mean   ${\bf 0 }$ and   direct/cross covariance  functions
   $$ C_{ij} (\vartheta; t) =  \exp (b_{ij} (t) \cos \vartheta) \cos ( b_{ij} (t)  \sin \vartheta),  ~~~~ \vartheta \in [0, \pi], ~ t \in \T,  ~ i, j =1, \ldots, m, $$
    which follows from    (\ref{cov.matrix.fun.2.0}) and the identity (see, e.g., (5) on page 98 of  \cite{Mangulis1965})
     $$ \sum_{n=0}^\infty \frac{b^n}{n!}  \cos (n \theta) =   \exp (b \cos \theta) \cos ( b \sin \theta),  
            ~~~~ b \in \R,  \theta \in [0, \pi]. $$

 \section{Concluding Remarks}

 Our focus is mostly on the spherical domain, although the vector random field in this paper has a spatio-temporal domain.
 While the temporal domain $\T$ is assumed to be either $\Z$ or $\R$, the results deduced here may  be easily extended to other cases. 
 
 The spatial domain $\S^d$ may be substtituted by  a $d$-dimensional compact two-point homogeneous Riemannian manifold  $\mathbb{M}^d$. It   is a compact Riemannian symmetric space of rank one, and  belongs to one of the
  following categories  (\cite{Helgason2011}, \cite{Wang1952}): the unit spheres $\S^d$ ($ d =1, 2, \ldots$), the real projective spaces $\mathbb{P}^d(\R)$ ($d = 2, \ldots$),
  the complex projective spaces  $\mathbb{P}^d(\mathbb{C})$ ($d = 4, 6, \ldots$), the quaternionic  projective spaces 
  $\mathbb{P}^d(\mathbb{H})$ ($d = 8, 12, \ldots$), and the Cayley  elliptic plane $\mathbb{P}^{16}  (Cay)$.  For the lowest dimensions,    $\mathbb{P}^1(\R) = \S^1, \mathbb{P}^2 (\mathbb{C}) = \S^2$,  and
  $\mathbb{P}^4(\mathbb{H}) = \S^4$.
  A  series representation of  a continuous and isotropic covariance function on $\mathbb{M}^d$  may be found  
  in \cite{Askey1976},   \cite{Gangolli1967}, \cite{Malyarenko2013}, with the ultraspherical polynomials   substituted by Jacobi polynomials, which  include the ultraspherical polynomial as a special case. 
 The general form  like those in Theorem 1 may be deduced for the covariance matrix function of an $m$-variate   random field  $\{  \bZ (\x; t),   \x \in \mathbb{M}^d,  t \in   \T \}$ that is isotropic on $\mathbb{M}^d$ and stationary on
 $\T$, although the approach in the proof of Theorem 1 may not be adopted, where Lemma 1 plays a key role.
 For an associated random field, is it possible to establish a series representation like (\ref{stoc1})
  with  $\S^d \times \T$  substituted by   $\mathbb{M}^d \times \T$?   This  would highly depend on
 whether an  orthogonal property   like that in Lemma 1 holds for  Jacobi polynomials  over $\mathbb{M}^d$ \cite{Ma2015}.

 Theorem 2 characterizes a  covariance matrix function  on $\S^\infty \times \T$,  whose entries are isotropic and continuous on $\S^\infty$ and stationary on $\T$. 
For an associated random field,  it would be of interest to   derive a series representation like (\ref{stoc1}). 
In the scalar and purely spherical case, a series representation of an  associated Gaussian random field  is  given by \cite{Berman1980},
        $$ Z (\x) = \sum_{n=0}^\infty b_n  Y_n (\x),  ~~~~~~~  \x  \in \S^\infty, $$
        where $\{ b_n, n \in \N_0 \}$ is a summable sequence of nonnegative numbers, for each $n \in \N_0$,   $\{ Y_n (\x),  \x  \in \S^\infty \}$ is a Gaussian random field with mean 0 and covariance
           $$ \cov ( Y_n (\x_1), Y_n (\x_2) ) = \cos^n (\vartheta (\x_1, \x_2)),   ~~~ \x_1, \x_2 \in  \S^\infty,  $$
           and $\{ Y_n (\x), \x \in \S^\infty \}$, $n \in \N_0$,   are mutually independent.
           Theorem 4 on page 83 of \cite{Yadrenko1983} gives an approach 
  to generate each $Y_n (\x)$  on $\S^\infty$, while  two generating methods   are available on $\S^d$, one   in         
   Subsection 12.3 of  \cite{Cohen2012} and the other  in Lemma 3.

   The ultraspherical functions $P_n^{ \left( \frac{d-1}{2} \right) } (\cos \v),  n \in \N_0$, are the basic spherical harmonics on $\S^d$
   ($d \ge 2$), analogous to $\cos (n \v)$ on $\S^1$.
   For every spherical harmonic $S_n (\x)$, it is possible to choose $h(n)$ points  $\y_1, \ldots, \y_{h(n)}$ on $\S^d$ such that
   $S_n (\x)$ is a linear combination of $P_n^{ \left( \frac{d-1}{2} \right) } (\x'\y_k), k =1, \ldots, h(n)$, according to  Theorem 9.6.4 of 
   \cite{Andrews1999}.  With such substitutions,  (\ref{Mokljacuk1979.exp}) can be rewritten in terms of the ultraspherical polynomials,
   although it is more completed than   (\ref{stoc1}) in the scalar case.
   At each level $n$,  only one term gets involved in (\ref{stoc1}), in contrast to $h(n)$ terms in   (\ref{Mokljacuk1979.exp}).
  Intuitively,  employing  finitely truncated expansions of (\ref{stoc1})  for approximation or simulation  would be more efficient than that of  (\ref{Mokljacuk1979.exp}). 
 An examination of   the convergent rate would be expected if finitely truncated expansions  of (\ref{stoc1}) are used for approximation or simulation.
   A purely spatial case with $d=2$ is studied in \cite{Lang2015}  with respect to the spectral representation (\ref{Mokljacuk1979.exp}).

 \section{Proofs}

\subsection{Proof of Theorem 1}

For a fixed $t \in \T$, consider two purely spatial random fields $\left\{  \bZ (\x; 0) + \bZ (\x; t),  \x \in \S^d \right\}$ and
$\left\{   \bZ (\x; 0) -  \bZ (\x; t),  \x \in \S^d \right\}$.  In terms of $\bC( \vartheta; t)$, their covariance matrix functions are, respectively,
     \begin{eqnarray*}
     &   &   \cov \left(     \bZ (\x_1; 0) + \bZ (\x_1; t), ~ \bZ (\x_2; 0) + \bZ (\x_2; t)       \right)     \\
     & = &  2  \bC (\vartheta (\x_1, \x_2); 0) +  \bC (\vartheta (\x_1, \x_2); t) +  \bC (\vartheta (\x_1, \x_2); -t),
     \end{eqnarray*}
     and
     \begin{eqnarray*}
     &   &   \cov \left(     \bZ (\x_1; 0) -  \bZ (\x_1; t), ~  \bZ (\x_2; 0) -  \bZ (\x_2; t)       \right)     \\
     & = &  2 \bC (\vartheta (\x_1, \x_2); 0) -  \bC (\vartheta (\x_1, \x_2); t) -  \bC (\vartheta (\x_1, \x_2); -t), 
                     ~~~~~ \x_1, \x_2 \in \S^d.
     \end{eqnarray*}
 We consider the case $d \ge 2$ only, while   a similar argument applies to the case $d=1$.     By Theorem 1 of  \cite{Ma2012}, these two covariance matrix functions must take the form
            \begin{equation}
            \label{eq1.thm1}
               2  \bC (\vartheta (\x_1, \x_2); 0) +  \bC (\vartheta (\x_1, \x_2); t) +  \bC (\vartheta (\x_1, \x_2); -t)  
                 =    \sum_{n=0}^\infty \mathbf{B}_{n+} (t)   P_n^{ \left( \frac{d-1}{2} \right) } (\cos \vartheta (\x_1, \x_2)),   
           \end{equation}
             \begin{equation}
            \label{eq2.thm1}
               2  \bC (\vartheta (\x_1, \x_2); 0) -  \bC (\vartheta (\x_1, \x_2); t) -  \bC (\vartheta (\x_1, \x_2); -t) 
              =  \sum_{n=0}^\infty \mathbf{B}_{n-} (t)   P_n^{ \left( \frac{d-1}{2} \right) } (\cos \vartheta (\x_1, \x_2)),
     \end{equation}
  where $ \mathbf{B}_{n+} (t) $ and  $ \mathbf{B}_{n-} (t) $  ($ n \in \N_0$)  are $m \times m$ positive definite matrices, and  
  $\sum\limits_{n=0}^\infty  \mathbf{B}_{n+}  (t)  P_n^{ \left( \frac{d-1}{2} \right) } (1)$  and  $\sum\limits_{n=0}^\infty  \mathbf{B}_{n-}  (t)  P_n^{ \left( \frac{d-1}{2} \right) } (1)$ converge.
  Taking the difference between (\ref{eq1.thm1}) and (\ref{eq2.thm1}) results in (\ref{cov.matrix.fun.1}), with 
        $$    \mathbf{B}_n (t) = \frac{1}{4}  \mathbf{B}_{n+} (t) - \frac{1}{4}  \mathbf{B}_{n-} (t),  ~~~~~~~~   n \in \N_0. $$
        Clear, $\mathbf{B}_n(t)$ is symmetric,
        and   $\sum\limits_{n=0}^\infty  \mathbf{B}_n (t)  P_n^{ \left( \frac{d-1}{2} \right) } (1)$ converges. 
        
        What remains is to verify that $ \mathbf{B}_n (t), t \in \T, $ is a  stationary covariance matrix function, for each fixed $n \in \N_0$.
        To this end, consider an $m$-variate  stochastic process
              $$ \mathbf{W}_n (t) = \int_{\S^d}  \frac{ \bZ (\x; t)+\tilde{\bZ} (\x; -t)}{\sqrt{2}}  P_n^{ \left( \frac{d-1}{2} \right) } (\x' \mathbf{U}) d \x,  ~~~~~~~~~~~~~~ t \in \T, $$
        where $\{ \tilde{\bZ} (\x; t), \x \in \S^d, t \in \T \}$ is an independent copy of  $\{ \bZ (\x; t), $  $ \x  \in \S^d, t \in \T \}$,  $\mathbf{U}$ is an $(d+1)$-variate random vector uniformly distributed on $\S^d$,  and $\mathbf{U}$,  $\{ \bZ (\x; t), $  $ \x  \in \S^d, t \in \T \}$ and  $\{ \tilde{\bZ} (\x; t), \x \in \S^d, t \in \T \}$ are independent.
         
         The  mean function of $\{ \mathbf{W}_n (t),  t \in \T \}$  is
          \begin{eqnarray*}
           \rE   \mathbf{W}_n (t)  &  = &   \rE  \int_{\S^d}  \frac{ \bZ (\x; t)+\tilde{\bZ} (\x; -t)}{\sqrt{2}}  P_n^{ \left( \frac{d-1}{2} \right) } (\x' \mathbf{U}) d \x    \\
             & =  &  \frac{1}{\omega_d} \int_{\S^d}   \int_{\S^d}  \rE \left(  \frac{ \bZ (\x; t)+\tilde{\bZ} (\x; -t)}{\sqrt{2}}  P_n^{ \left( \frac{d-1}{2} \right) } (\x' \mathbf{u})  \right) d \x   d \mathbf{u}  \\
            & =  &  \frac{\sqrt{2}  \rE \bZ (\x; t)  }{\omega_d} \int_{\S^d}   \int_{\S^d}    P_n^{ \left( \frac{d-1}{2} \right) } (\x' \mathbf{u})   d \x   d \mathbf{u}  \\
             & =  &  \left\{
                          \begin{array}{ll}
                          \sqrt{2} \omega_d  \rE \bZ (\x; t),   ~   &    ~  n= 0,   \\
                          0,     ~   &     ~  n \in \N,
                          \end{array}    \right. 
                  \end{eqnarray*}
              where the last equality follows from Lemma 1.     
            As is shown above,    the covariance matrix function of  an $m$-variate random field $\left\{  \frac{ \bZ (\x; t)+\tilde{\bZ} (\x; -t)}{\sqrt{2}}, \x \in \S^d, t \in \T \right\}$  is of the form
                   \begin{eqnarray*}
                    &     &  \cov \left( \frac{ \bZ (\x_1; t_1)+\tilde{\bZ} (\x_1; -t_1)}{\sqrt{2}}, ~  \frac{ \bZ (\x_2; t_2)+\tilde{\bZ} (\x_2; -t_2)}{\sqrt{2}}  \right)   \\
                    & = &   \frac{ \bC (\vartheta (\x_1, \x_2); t_1-t_2)+\bC ( \vartheta (\x_1, \x_2); t_2-t_1)}{2}   \\
                    & = &   \sum_{k=0}^\infty \mathbf{B}_{k} (t_1-t_2)   P_k^{ \left( \frac{d-1}{2} \right) } (\cos \vartheta (\x_1, \x_2))   \\
                    & = &   \sum_{k=0}^\infty \mathbf{B}_{k} (t_1-t_2)   P_k^{ \left( \frac{d-1}{2} \right) } (\x'_1 \x_2),     ~~~~~~~ \x_1, \x_2 \in \S^d, t_1, t_2 \in \T. 
                    \end{eqnarray*}
           From this observation and Lemma 1 we obtain  the covariance matrix function of $\{ \mathbf{W}_n (t), t \in \T \}$,
                     \begin{eqnarray*}
                   &   &   \cov ( \mathbf{W}_n(t_1), ~  \mathbf{W}_n (t_2) )  \\
                   & = &  \cov   \left(       \int_{\S^d}  \frac{ \bZ (\x; t_1)+\tilde{\bZ} (\x; -t_1)}{\sqrt{2}}  P_n^{ \left( \frac{d-1}{2} \right) } (\x' \mathbf{U}) d \x, ~~  
                                           \int_{\S^d}  \frac{ \bZ (\y; t_2)+\tilde{\bZ} (\y; -t_2)}{\sqrt{2}}  P_n^{ \left( \frac{d-1}{2} \right) } (\y' \mathbf{U}) d \y                           \right) \\
                    & = &  \frac{1}{\omega_d}   \int_{\S^d}  \cov   \left(       \int_{\S^d}  \frac{ \bZ (\x; t_1)+\tilde{\bZ} (\x; -t_1)}{\sqrt{2}}  P_n^{ \left( \frac{d-1}{2} \right) } (\x' \mathbf{u}) d \x, ~~
                                           \int_{\S^d}  \frac{ \bZ (\y; t_2)+\tilde{\bZ} (\y; -t_2)}{\sqrt{2}}  P_n^{ \left( \frac{d-1}{2} \right) } (\y' \mathbf{u}) d \y                           \right)   d \mathbf{u}  \\
                    & = &  \frac{1}{\omega_d}   \int_{\S^d}      \int_{\S^d}        \int_{\S^d}   \cov \left( \frac{ \bZ (\x; t_1)+\tilde{\bZ} (\x; -t_1)}{\sqrt{2}}, ~  \frac{ \bZ (\y; t_2)+\tilde{\bZ} (\y; -t_2)}{\sqrt{2}}  \right)
                       P_n^{ \left( \frac{d-1}{2} \right) } (\x' \mathbf{u})    P_n^{ \left( \frac{d-1}{2} \right) } (\y' \mathbf{u})  d \x d \y                             d \mathbf{u}  \\
                    & = &     \frac{1}{\omega_d}   \int_{\S^d}         \int_{\S^d}     \int_{\S^d}    \frac{ \bC (\vartheta (\x, \y); t_1-t_2)+\bC ( \vartheta (\x, \y); t_2-t_1)}{2}
                       P_n^{ \left( \frac{d-1}{2} \right) } (\x' \mathbf{u})    P_n^{ \left( \frac{d-1}{2} \right) } (\y' \mathbf{u})  d \x d \y                             d \mathbf{u}  \\
                     & = &     \frac{1}{\omega_d}   \int_{\S^d}        \int_{\S^d}      \int_{\S^d}      \sum_{k=0}^\infty \mathbf{B}_{k} (t_1-t_2)   P_k^{ \left( \frac{d-1}{2} \right) } (\x' \y)
                        P_n^{ \left( \frac{d-1}{2} \right) } (\x' \mathbf{u})    P_n^{ \left( \frac{d-1}{2} \right) } (\y' \mathbf{u})  d \x d \y                             d \mathbf{u}  \\
                    & = &     \frac{1}{\omega_d}     \sum_{k=0}^\infty \mathbf{B}_{k} (t_1-t_2)      \int_{\S^d}   \left\{     \int_{\S^d}    \left(     \int_{\S^d}  
                         P_k^{ \left( \frac{d-1}{2} \right) } (\x' \y)
                        P_n^{ \left( \frac{d-1}{2} \right) } (\x' \mathbf{u})   d \x \right)   P_n^{ \left( \frac{d-1}{2} \right) } (\y' \mathbf{u})   d \y               \right\}               d \mathbf{u}  \\
                    & = &     \frac{1}{\omega_d}     \mathbf{B}_{n} (t_1-t_2)      \int_{\S^d}   \left(   \frac{\omega_d}{\alpha_n^2}  \int_{\S^d}      
                         P_n^{ \left( \frac{d-1}{2} \right) } (\y' \mathbf{u})   P_n^{ \left( \frac{d-1}{2} \right) } (\y' \mathbf{u})   d \y               \right)              d \mathbf{u}  \\
                   & = &     \frac{1}{\omega_d}     \mathbf{B}_{n} (t_1-t_2)      \int_{\S^d}   \left(  \frac{\omega_d}{\alpha_n^2} \right)^2      
                         P_n^{ \left( \frac{d-1}{2} \right) } (1)                d \mathbf{u}  \\
                    & = &        \mathbf{B}_{n} (t_1-t_2)       \left(  \frac{\omega_d}{\alpha_n^2} \right)^2      
                         P_n^{ \left( \frac{d-1}{2} \right) } (1),    ~~~~~~ t_1, t_2 \in \T, 
                                             \end{eqnarray*}
             which implies that  $ \mathbf{B}_n (t)$ is a  stationary covariance matrix function on $\T$.

               \subsection{Proof of Theorem 2}

                 (i)     Suppose that (\ref{cov.matrix.fun.2.0}) is the covariance matrix function of an $m$-variate random field $\{ \bZ (\x; t), \x \in \S^d, t \in \T \}$.
               The existence of $\bC (0; 0)$ ensures   the convergence of   $\sum\limits_{n=0}^\infty  \mathbf{B}_n (0)$. 
               To verify that $\mathbf{B}_n(t)$ is a stationary covariance matrix function on $\T$ for
               each fixed $n \in \N_0$,  consider an $m$-variate   stochastic process
                  $$ \mathbf{W}_n (t) =  \frac{1}{\pi} \int_0^{2 \pi}  \bZ (\x; t)  \cos (n \theta)  d \theta,  ~~~~~~~~~ t \in \T, $$
              where    $\x = ( \cos \theta, \sin \theta)' \in \S^1$,   $0 \le \theta \le 2 \pi$.
             The covariance matrix function of $\{ \mathbf{W}_n (t), t \in \T \}$ is given by 
                   \begin{eqnarray*}
                  &   &    \cov ( \mathbf{W}_n (t_1),   \mathbf{W}_n (t_2) )    \\
                   & =  &   \frac{1}{\pi^2}   \cov \left(     \int_0^{2 \pi}  \bZ (\x_1; t_1)  \cos (n \theta_1) d \theta_1, 
                               ~ 
                                 \int_0^{2 \pi}  \bZ (\x_2; t_2)  \cos (n \theta_2)  ) d \theta_2    \right)    \\
                 & =  &   \frac{1}{\pi^2}      \int_0^{2 \pi}   \int_0^{2 \pi}     \cov ( \bZ (\x_1; t_1),  \bZ (\x_2; t_2) )   \cos (n \theta_1) \cos (n \theta_2)     d \theta_1  d \theta_2       \\
                & = &       \frac{1}{\pi^2}      \int_0^{2 \pi}   \int_0^{2 \pi}  \sum\limits_{k=0}^\infty \mathbf{B}_k (t_1-t_2) \cos (k \vartheta (\x_1,  \x_2))  \cos (n \theta_1) \cos (n \theta_2) 
                         d \theta_1  d \theta_2       \\
                   & = &       \frac{1}{\pi^2}  \sum\limits_{k=0}^\infty \mathbf{B}_k (t_1-t_2)     \int_0^{2 \pi}   \int_0^{2 \pi}   \cos (k (\theta_1-\theta_2))  \cos (n \theta_1) \cos (n \theta_2) 
                         d \theta_1  d \theta_2       \\
                               & = &      \frac{1}{\pi^2}   \mathbf{B}_n (t_1-t_2)     \int_0^{2 \pi}   \int_0^{2 \pi}   \cos (n (\theta_1-\theta_2))  \cos (n \theta_1) \cos (n \theta_2) 
                         d \theta_1  d \theta_2       \\
                              & =  &  \left\{
                                          \begin{array}{ll}
                                          2 \mathbf{B}_0 (t_1-t_2),    ~   &   ~ n = 0,    \\
                                          \mathbf{B}_n (t_1-t_2),  ~   &    ~ n \in \N,  ~~~~~~ t_1, t_2 \in \T,
                                          \end{array}    \right. 
                                            \end{eqnarray*}
            which implies that  $ \mathbf{B}_n (t)$ is a  stationary covariance matrix function on $\T$.

          Conversely,  if   $ \mathbf{B}_n (t)$  ($n \in \N_0$) are  stationary covariance matrix functions on $\T$ and $\sum\limits_{n=0}^\infty  \mathbf{B}_n (0)$ converges, then, as Theorem 5 shows, we can generate
          an $m$-variate random field with  (\ref{cov.matrix.fun.2.0}) as its covariance matrix function, so that  (\ref{cov.matrix.fun.2.0}) satisfies  inequality (\ref{positive.definite}). By Theorem 8 of \cite{Ma2011},
          there exists an  $m$-variate  Gaussian or elliptically contoured random field with  (\ref{cov.matrix.fun.2.0}) as its covariance matrix function.

           (ii)    We give a proof of the ``only if" part here, while the ``if " part is analogous to that in the proof of Part (i).   Suppose that (\ref{cov.matrix.fun.2.1}) is the covariance matrix function of an $m$-variate random field $\{ \bZ (\x; t), \x \in \S^d, t \in \T \}$.
               Evidently, the existence of $\bC (0; 0)$ implies the convergence of   $\sum\limits_{n=0}^\infty  \mathbf{B}_n (0)  P_n^{ \left( \frac{d-1}{2} \right) } (1)$. 
              For each fixed $n \in \N_0$,  consider an $m$-variate  stochastic process
                  $$ \mathbf{W}_n (t) = \int_{\S^d} \bZ (\x; t)  P_n^{ \left( \frac{d-1}{2} \right) } (\x' \mathbf{U})  d \x,  ~~~~~~~~~ t \in \T, $$
              where $\mathbf{U}$ is a ($d+1$)-dimensional random vector uniformly distributed on $\S^d$ and independent with     $\{ \bZ (\x; t), \x \in \S^d, t \in \T \}$.
              In a way similar to  the proof of Theorem 1, we apply Lemma 1 to  obtain that the covariance matrix function of $\{ \mathbf{W}_n (t), t \in \T \}$ is positively propositional to 
              $ \mathbf{B}_n (t)$. More precisely,
                   \begin{eqnarray*}
                    \cov ( \mathbf{W}_n (t_1),   \mathbf{W}_n (t_2) )  
                   & =  &    \mathbf{B}_{n} (t_1-t_2)       \left(  \frac{\omega_d}{\alpha_n^2} \right)^2      
                         P_n^{ \left( \frac{d-1}{2} \right) } (1),    ~~~~~~ t_1, t_2 \in \T, 
                   \end{eqnarray*}
            so that  $ \mathbf{B}_n (t)$ is a  stationary covariance matrix function on $\T$.

               \subsection{Proof of Theorem 3}

For a fixed $t \in \T$,  in a way similar to that in the proof of Theorem 1 it can be verify that 
     $  2 \bC (\vartheta (\x_1, \x_2); 0) +  \bC (\vartheta (\x_1, \x_2); t) +  \bC (\vartheta (\x_1, \x_2); -t)$
     and  $ 2 \bC (\vartheta (\x_1, \x_2); 0) -  \bC (\vartheta (\x_1, \x_2); t) - \bC (\vartheta (\x_1, \x_2); -t)$
     are isotropic covariance matrix functions on $\S^\infty$.  They  necessarily take the form,     by Theorem 4 of  \cite{Ma2015}, 
            \begin{equation}
            \label{eq1.thm2}
                 2 \bC (\vartheta (\x_1, \x_2); 0) +  \bC (\vartheta (\x_1, \x_2); t) +  \bC (\vartheta (\x_1, \x_2); -t)  
                 =    \sum_{n=0}^\infty \mathbf{B}_{n+} (t)   \cos^n ( \vartheta (\x_1, \x_2)),   
           \end{equation}
             \begin{equation}
            \label{eq2.thm2}
                2 \bC (\vartheta (\x_1, \x_2); 0) -  \bC (\vartheta (\x_1, \x_2); t) -  \bC (\vartheta (\x_1, \x_2); -t) 
              =  \sum_{n=0}^\infty \mathbf{B}_{n-} (t)   \cos^n  (\vartheta (\x_1, \x_2)),
     \end{equation}
  where $ \mathbf{B}_{n+} (t) $ and  $ \mathbf{B}_{n-} (t) $  ($ n \in \N_0$)  are $m \times m$ positive definite matrices, and  
  $\sum\limits_{n=0}^\infty  \mathbf{B}_{n+}  (t) $  and  $\sum\limits_{n=0}^\infty  \mathbf{B}_{n-}  (t)  $ converge.
  The representation (\ref{cov.matrix.fun.3}) results from  taking  the difference between (\ref{eq1.thm2}) and (\ref{eq2.thm2}), and 
        $$    \mathbf{B}_n (t) =  \frac{1}{4}  \mathbf{B}_{n+} (t) - \frac{1}{4}  \mathbf{B}_{n-} (t),  ~~~~~~~~   n \in \N_0. $$
        Clear, $\mathbf{B}_n(t)$ is symmetric,
        and   $\sum\limits_{n=0}^\infty  \mathbf{B}_n (t)  $ converges.

        In particular,  $\mathbf{B}_0 (t) = \bC (0; t), t \in \T, $ is a stationary covariance matrix function.  For each $n \in \N$,  we are going to confirm  that $ \mathbf{B}_n (t), t  \in \T, $ is a  stationary covariance matrix function.  For every $d \ge 2$,
        a version (\ref{cov.matrix.fun.2.1}) of $\bC (\vartheta; t)$ is derived from   (\ref{cov.matrix.fun.3})  by using the formula    (\ref{Binghamlemma1}), 
            \begin{eqnarray*}
            \bC (\vartheta; t) & = &  \sum_{n=0}^\infty \mathbf{B}_n  (t)   \cos^n ( \vartheta (\x_1, \x_2))   \\
                 & = &    \sum_{n=0}^\infty \mathbf{B}_n  (t)   \sum\limits_{k=0}^{ \left[ \frac{n}{2} \right]}
                   \beta_{k,n}^{\left( \frac{d-1}{2} \right)} P_{n-2 k}^{ \left( \frac{d-1}{2} \right) } ( \cos \vartheta )   \\
                  & = &    \sum_{n=0}^\infty  \mathbf{A}_n^{ \left( \frac{d-1}{2} \right) }   (t)     P_{n}^{ \left( \frac{d-1}{2} \right) } ( \cos \vartheta ),  ~~~~~~ \vartheta \in [0, \pi], ~ t \in \T,    
            \end{eqnarray*}
            where
            $$      \mathbf{A}_n^{ \left( \frac{d-1}{2} \right) }   (t)  
                            = \sum\limits_{k=0}^\infty \beta_{k, 2 k+n}^{\left( \frac{d-1}{2} \right)}  \mathbf{B}_{2k+n}   (t),  ~~~~~~  n \in \N_0. $$  
                            Since  a covariance matrix function $ \bC (\vartheta; t)$ on $\S^\infty \times \T$ is also a covariance matrix function on $\S^d \times \T$ for every $d \in \N$,
                            applying Theorem 1 to  $ \bC (\vartheta; t)$ on  $\S^d \times \T$ we obtain that, for each $n \in \N_0$, $\mathbf{A}_n^{ \left( \frac{d-1}{2} \right) }   (t) $ is a stationary
                            covariance matrix function on $\T$.    So is  $ \frac{\mathbf{A}_n^{ \left( \frac{d-1}{2} \right) }   (t) }{ \beta_{0, n}^{\left( \frac{d-1}{2} \right)} },$ $ t \in \T$, by Theorem 6 of \cite{Ma2011}.
                            
                          For $k \in \N$,   it follows from the formula $\Gamma (x+1) = x \Gamma (x)$ that 
                               \begin{eqnarray*}
                                        \frac{ \beta_{k, 2k+n}^{\left( \frac{d-1}{2} \right)}}{\beta_{0,n}^{\left( \frac{d-1}{2} \right)} }    
                                   & = &  \frac{(n+2k)!}{2^k n! k!} \frac{n+k+\frac{d-1}{2}}{n+\frac{d-1}{2} }
                                             \frac{ \Gamma \left( n+\frac{d+1}{2} \right)}{ \Gamma \left( n+k +\frac{d+1}{2}  \right)}   \\
                                   & =  &    \frac{(n+2k)!}{2^k n! k!} \frac{n+k+\frac{d-1}{2}}{n+\frac{d-1}{2} }  \prod_{l=0}^k  \left( n+l+\frac{d+1}{2}  \right)^{-1}    \\
                                   &  \to &  0,    ~~~~~~~~~~ \mbox{as} ~~ d \to \infty.   
                                          \end{eqnarray*}
                            Hence,
                            $ \lim\limits_{d \to \infty}   \frac{ \mathbf{A}_n^{ \left( \frac{d-1}{2} \right) }   (t) } { \beta^{\left( \frac{d-1}{2} \right)} _{0, n } }    =  \mathbf{B}_n (t)    $
                            is a  stationary
                            covariance matrix function on $\T$.   
                            
             (ii)    The ``only if" part follows from Part (i), and  the ``if" part from    Theorem 8 of \cite{Ma2011}.

  \subsection{Proof of Theorem 4}
  
   The convergent  assumption  of $\sum\limits_{n=0}^\infty {\bf B}_n (0)  P_n^{ \left( \frac{d-1}{2} \right) } (1)$
       ensures  the mean square convergence of the series at the right hand of (\ref{stoc1}).   
                        In fact,  for $n_1, n_2 \in \N$, we have
           \begin{eqnarray*}
      &  &      \rE \left( \sum_{i=n_1}^{n_1+n_2}  \bV_i (t)  P_i^{ \left( \frac{d-1}{2} \right) } (\x' {\bf U}) \right)
                  \left( \sum_{j=n_1}^{n_1+n_2}  \bV_j (t) P_j^{ \left( \frac{d-1}{2} \right) } (\x' {\bf U}) \right)' \\
      & =  &      \rE \left(  \sum_{i=n_1}^{n_1+n_2}  \sum_{j=n_1}^{n_1+n_2} 
           \bV_i (t)   \bV'_j  (t)   P_i^{ \left( \frac{d-1}{2} \right) } (\x' {\bf U})  P_j^{ \left( \frac{d-1}{2} \right) } (\x' {\bf U}) \right)' \\    
      & =  &       \sum_{i=n_1}^{n_1+n_2}  \sum_{j=n_1}^{n_1+n_2} 
           \rE ( \bV_i  (t)  \bV'_j (t) )   \rE \left( P_i^{ \left( \frac{d-1}{2} \right) } (\x' {\bf U})  P_j^{ \left( \frac{d-1}{2} \right) } (\x' {\bf U})   \right) \\
      & =  &    \omega_d    \sum_{i=n_1}^{n_1+n_2} 
           {\bf B}_i (0)    P_i^{ \left( \frac{d-1}{2} \right) } (1)    \\
      &   \to  &  0,     ~~~~~~~  \mbox{as} ~  n_1, n_2 \to \infty,     
            \end{eqnarray*}
              where the second equality follows from   the independent assumption between ${\bf U}$ and $\{ \bV_n (t),  t \in \T \}$,  and  the third  one from Lemma 2.  
           
 Under the independent assumption among ${\bf U}$ and $\{ \bV_n (t), t \in \T \}$, $n \in \N_0$, we obtain the mean and covariance matrix functions of
  $\{ \bZ (\x; t), \x \in \S^d, t \in \T \}$ from Lemma 2, with
     $$ \rE \bZ (\x; t ) = \sum_{n=0}^\infty     \rE \bV_n  (t) \rE P_n^{ \left( \frac{d-1}{2} \right) } (\x' {\bf U})  = {\bf 0},
        ~~~ \x \in \S^d,  t \in \T, $$
   and
      \begin{eqnarray*}
      &      &   \cov ( \bZ (\x_1; t_1), \bZ ( \x_2; t_2) ) \\
      &  = &   \cov \left(   \sum_{i=0}^\infty   \bV_i (t_1)   P_i^{ \left( \frac{d-1}{2} \right) } (\x' {\bf U}),
                                     ~ \sum_{j=0}^\infty   \bV_j (t_2)  P_j^{ \left( \frac{d-1}{2} \right) } (\x' {\bf U})  \right)  \\
      &  = &   \sum_{i=0}^\infty  \sum_{j=0}^\infty    \rE ( \bV_i (t_1) \bV'_j (t_2) ) 
                        \rE \left( 
                                           P_i^{ \left( \frac{d-1}{2} \right) } (\x'_1 {\bf U}) P_j^{ \left( \frac{d-1}{2} \right) } (\x'_2 {\bf U}) \right)   \\
      &  = &   \sum_{n=0}^\infty    {\bf B}_n  (t_1-t_2)    \cov \left( 
                                         \alpha_n  P_n^{ \left( \frac{d-1}{2} \right) } (\x'_1 {\bf U}), ~ \alpha_n P_n^{ \left( \frac{d-1}{2} \right) } (\x'_2 {\bf U}) \right)   \\                                     
     &  = &   \sum_{n=0}^\infty     {\bf B}_n  (t_1-t_2) 
                                           P_n^{ \left( \frac{d-1}{2} \right) } ( \cos \vartheta (\x_1,  \x_2) ),  ~~~~~~~~~~ \x_1, \x_2 \in \S^d, ~ t_1, t_2 \in \T.
                                            \end{eqnarray*}     
                    The latter is obviously isotropic and continuous on $\S^d$ and stationary on $\T$.

  \subsection{Proof of Theorem 5}
  
    The series at the right hand side of (\ref{stoc2}) is convergent in mean square,  since $\sum\limits_{n=0}^\infty \mathbf{B}_n (0)$ is convergent and,  for $n_1, n_2 \in \N_0$,
              \begin{eqnarray*}
      &  &      \rE \left( \sum_{i=n_1}^{n_1+n_2}  ( \bV_{i1} (t)  \cos (i \theta)  + \bV_{i2}  (t) \sin (i \theta)   )
                 \right)
                  \left( \sum_{j=n_1}^{n_1+n_2}  ( \bV_{j1} (t)  \cos (j \theta)  + \bV_{j2}  (t) \sin (j \theta)   ) \right)' \\
      & =  &     \sum_{i=n_1}^{n_1+n_2}  \sum_{j=n_1}^{n_1+n_2}  \rE \left\{ 
            \bV_{i1}(t)   \bV'_{j1} (t)    \cos (i \theta) \cos (j \theta) 
        +   \bV_{i1}  (t)  \bV'_{j2}  (t)    \cos (i \theta) \sin (j \theta)  \right. \\    
      &   &       \left.       + \bV_{i2} (t)   \bV'_{j1}  (t)   \sin (i \theta) \cos (j \theta)  
            +  \bV_{i2}  (t)  \bV'_{j2}  (t)   \sin (i \theta) \sin (j \theta)  \right\}  \\    
      & =  &     \sum_{n=n_1}^{n_1+n_2}    
           {\bf B}_n  (0) \ \\    
               &  \to  &   {\bf 0},                  ~~~~~\mbox{as} ~~~ n_1 \to \infty, ~ n_2 \to \infty,  
            \end{eqnarray*}
              where the second equality is due to  the assumptions on $ \{ {\bf V}_{n1} (t),  t \in \T \}$  and $\{ {\bf V}_{n2} (t), t \in \T \}$,  $ n \in \N_0 \}$.
              
 Clearly, the mean  function of
  $\{ \bZ (\x; t), \x \in \S^1, t \in \T \}$ is identical to ${\bf 0}$, and its covariance matrix function is
    \begin{eqnarray*}
      &      &   \cov ( \bZ (\x_1; t_1), \bZ ( \x_2; t_2) ) \\
      &  = &   \cov \left(   \sum_{i=0}^\infty     ( \bV_{i1} (t_1)  \cos (i \theta_1)  + \bV_{i2} (t_1)  \sin (i \theta_1)   ),
                                     ~   \sum_{j=0}^\infty     ( \bV_{j1} (t_2)  \cos (j \theta_2)  + \bV_{j2}  (t_2)  \sin (j\theta_2)   )
                                     \right)  \\
     &  = &   \sum_{i=0}^\infty   \sum_{j=0}^\infty \left\{  \cov ( \bV_{i1} (t_1),  \bV_{j1} (t_2) )  \cos (i \theta_1)  \cos (j \theta_2)
                          +  \cov (  \bV_{i2} (t_1),  \bV_{j2}  (t_2) )   \sin (i \theta_1)  \sin (j\theta_2) \right\}    \\
      &  = &   \sum_{n=0}^\infty  {\bf B}_n (t_1-t_2) \{ \cos (n \theta_1) \cos (n \theta_2)
                          + \sin  (n \theta_1) \sin (n \theta_2)   \}  \\
                          &  = &   \sum_{n=0}^\infty    {\bf B}_n  (t_1-t_2)     \cos (n (\theta_1 -\theta_2)) \\                                     
     &  = &   \sum_{n=0}^\infty     {\bf B}_n  
                                           \cos (n \vartheta (\x_1,  \x_2) ),  ~~~~~~~~~~ \x_k = (\cos \theta_k, \sin \theta_k)'
                                            \in \S^1,  ~ t_k \in \T,   ~~~~~ k = 1, 2.
                                            \end{eqnarray*}     
                                                                 
                      \section*{Acknowledgment}
 
    A reviewer's valuable comments and  helpful    suggestions are gratefully acknowledged.

\end{document}